\documentclass[12pt]{article}
\usepackage[utf8]{inputenc}
\usepackage[T1]{fontenc}
\usepackage[frenchb,english]{babel} 
\usepackage{textcomp}
\usepackage{amsmath,amssymb}
\usepackage{amsthm}
\usepackage{lmodern}

\usepackage[a4paper,margin=2cm]{geometry}

\usepackage{graphicx}             
\usepackage{xcolor}               
\usepackage{microtype,stmaryrd}          

\usepackage{hyperref}
\hypersetup{pdfstartview=XXZ}

\usepackage{bbm}
\usepackage{mathrsfs}
\usepackage{subfig}

\def\build#1_#2^#3{\mathrel{
\mathop{\kern 0pt#1}\limits_{#2}^{#3}}}

\newtheorem{theorem}{Theorem}
\newtheorem{proposition}[theorem]{Proposition}

\newtheorem{lemma}[theorem]{Lemma}

\def\w{\mathrm{w}}

\def\f{\mathcal{F}}
\def\r{\mathcal{R}}
\def\g{\mathcal{G}}
\def\W{\mathcal{W}}
\def\S{\mathcal{S}}

\def\N{\mathbb{N}}

\def\P{\mathbb{P}}

\def\E{\mathbb{E}}

\def\R{\mathbb{R}}

\def\n{\mathcal{N}}

\def\ve{{\varepsilon}}
\def\la{\longrightarrow}

\def\ov{\overline}

\def\dd{\mathrm{d}}
\def\wh{\widehat}
\def\wt{\widetilde}

\def\BP{\mathbb{B}\mathbb{P}}

\def\XX{\mathbf{X}}
\def\II{\mathrm{i}}

\def\rem{\noindent{\bf Remark. }}

\author{Jean-Fran\c cois Le Gall\footnote{e-mail: jean-francois.le-gall@universite-paris-saclay.fr}} 
\title{The area of spheres in the Brownian plane}
\date{\small Universit\'e Paris-Saclay}

\begin{document}
\maketitle

\abstract{We consider the area of spheres centered at the
distinguished point in the Brownian plane. As a function of the radius, the resulting process 
has continuously differentiable sample paths. Furthermore, the pair consisting of the process
and its derivative is time-homogeneous Markov and satisfies an explicit stochastic
differential equation.}

\section{Introduction}
\label{intro}

The Brownian plane $(\BP,D_\infty)$ is a random locally compact metric space, which is homeomorphic to the plane,
and appears as the scaling limit of various classes of random planar maps. In particular, the Brownian 
plane has been shown \cite{Bud,plane} to be the scaling limit of the infinite random lattices called the uniform infinite planar triangulation (UIPT) and the 
uniform infinite planar quadrangulation (UIPQ), which 
have been extensily studied. The Brownian plane comes with a distinguished point $\rho$ and 
a volume measure denoted by $\mathrm{Vol}(\cdot)$, which coincides with the Hausdorff measure
associated with an appropriate gauge function (this follows from the results of \cite{Hausdorff} and the coupling between
the Brownian plane and the Brownian sphere). 

In this work, we are interested in the 
profile of distances from $\rho$, which is the $\sigma$-finite measure 
$\Gamma$ on $\R_+$ defined by setting
\begin{equation}
\label{profi}
\Gamma(A):=\int_{\BP} \mathbf{1}_A(D_\infty(\rho,a))\,\mathrm{Vol}(\dd a),
\end{equation}
for every Borel subset $A$ of $\R_+$. 

\begin{proposition}
\label{localtimesBP}
The measure $\Gamma$
has a continuously differentiable density with respect to Lebesgue measure on $\R_+$. We denote this density by
$(L^\infty_t)_{t\in\R_+}$ and its derivative by $(\dot L^\infty_t)_{t\in\R_+}$. Moreover we have $L^\infty_0=\dot L^\infty_0=0$,  and  $L^\infty_t>0$ for every $t>0$, a.s.
\end{proposition}

Since, for every $t>0$,
$$L^\infty_t=\lim_{\ve \to 0} \frac{1}{\ve} \,\mathrm{Vol}(\{a\in\BP:t\leq D_\infty(\rho,a)\leq t+\ve\})$$
it is natural to interpret $L^\infty_t$ as the area of the sphere of radius $t$ centered at $\rho$
in $\BP$. It is essentially obvious that the process $(L^\infty_t)_{t\geq 0}$ is not Markovian:
for $t>0$, the derivative $\dot L^\infty_t$ is a function of the past $(L^\infty_s)_{0\leq s\leq t}$ that
gives more information on the future $(L^\infty_s)_{s\geq t}$ than the sole knowledge of the present $L^\infty_t$. 

Set $\f^\circ_t=\sigma(L^\infty_s:s\leq t)$  for every $t\geq 0$, and let $(\f_t)_{t\geq 0}$ be the completion of 
the filtration $(\f^\circ_t)_{t\geq 0}$.

\begin{theorem}
\label{Markov-LT}
The process $(L^\infty_s,\dot L^\infty_s)_{s\geq 0}$ is a time-homogeneous Markov process
with respect to the filtration $(\f_t)_{t\geq 0}$.
\end{theorem}

The other main contribution of the present work is to show that this Markov process satisfies an explicit stochastic differential equation. Before stating this result, we need 
to introduce some notation. For every $t>0$, we let $p_t$ be the (continuous) density at time $t$ of the stable L\'evy process with index $3/2$
and no negative jumps, whose law is characterized by its Laplace exponent $\psi(\lambda)=\sqrt{2/3}\,\lambda^{3/2}$. Alternatively, we
can characterize $p_t$ by its Fourier transform
$$\int_\R e^{\II u x}\,p_t(x)\,\dd x = \exp(-c_0t\, |u|^{3/2}\,(1+\II\,\mathrm{sgn}(u))),$$
where $c_0=1/\sqrt{3}$ and $\mathrm{sgn}(u)=\mathbf{1}_{\{u>0\}}-\mathbf{1}_{\{u<0\}}$. Then $x\mapsto p_t(x)$ is strictly positive, infinitely differentiable and 
has bounded derivatives. We write $p'_t$ for the derivative of this function.

Recall that 
a linear Brownian motion $B$ is an $(\f_t)$-Brownian motion if it is adapted to $(\f_t)_{t\geq 0}$ and has 
independent increments with respect to the filtration $(\f_t)_{t\geq 0}$.

\begin{theorem}
\label{mainT}
Set, for every $t>0$ and $x\in\R$,
$$h(t,x)= -8t\,\frac{p'_t(-x/2)}{p_t(-x/2)}+ \frac{4}{3}\,\frac{x^2}{t}.$$
Then for every $t>0$ and $x\in\R$, we have $h(t,x)>0$ and
$$\int_0^t h(L^\infty_s,\dot L^\infty_s)\,\dd s <\infty,\quad a.s.$$
There exists an $(\f_t)$-Brownian motion $(B_t)_{t\geq0}$ such that $(\dot L^\infty_t)_{t\geq 0}$
is a semimartingale satisfying the equation
\begin{equation}
\label{EQ1}
\dot L^\infty_t = 4\int_0^t\sqrt{L^\infty_s}\,\dd B_s + \int_0^t h(L^\infty_s,\dot L^\infty_s)\,\dd s.
\end{equation}
\end{theorem}

Since we have also trivially
\begin{equation}
\label{EQ2}
L^\infty_t=\int_0^t\dot L^\infty_s\,\dd s
\end{equation}
we can view \eqref{EQ1} and \eqref{EQ2} as a system of stochastic differential equations satisfied by the pair $(L^\infty_s,\dot L^\infty_s)_{s\geq 0}$.
One may obviously ask for uniqueness of the solution of this system. If $\ve>0$ is fixed, we know by Proposition \ref{localtimesBP} that $L^\infty_\ve>0$ a.s.
and then the classical uniqueness results for stochastic differential equations under Lipschitz conditions show that the solution starting at time $\ve$ from $(L^\infty_\ve,\dot L^\infty_\ve)$
is unique up to time $\inf\{t\geq\ve:L^\infty_t=0\}$. However the latter hitting time of $0$ is infinite a.s.~(by Proposition \ref{localtimesBP} again) and so
we obtain that $(L^\infty_t)_{t\geq \ve}$ is the unique solution of the system \eqref{EQ1} and \eqref{EQ2} that starts from $(L^\infty_\ve,\dot L^\infty_\ve)$ at time $\ve$. 
On the other hand, this argument does not rule out the possibility that there may exist other solutions starting from $(0,0)$ at time $0$ (see 
\cite{GLMNS} for examples of very similar systems where uniqueness starting from $(0,0)$ may fail). 

\smallskip
As discussed in the introduction of \cite{LGP}, the functions $p_t$ and $p'_t$ have an explicit expression in terms of the
classical Airy function $\mathrm{Ai}$. In particular,
\begin{equation} 
\label{formupt}
p_t(x)= 6^{-1/3}\,t^{-2/3} \,\mathcal{A}(6^{-1/3}t^{-2/3}x),
\end{equation}
where
\begin{equation}
\label{map-Airy}
\mathcal{A}(x)=-2\,e^{2x^3/3}\Big( x\mathrm{Ai}(x^2) + \mathrm{Ai}'(x^2)\Big).
\end{equation}
The function $x\mapsto \mathcal{A}(-x)$ is the density of the so-called map-Airy distribution (see \cite[Definition 1]{BFSS}). 
From these formulas (and the Airy equation $\mathrm{Ai}''(x)=x\,\mathrm{Ai}(x)$), one derives the following
expression for the drift $h$ in \eqref{EQ1}. We have 
\begin{equation}
\label{formu-drift}
h(t,x)=- 8\times6^{-1/3}t^{1/3}\, \kappa\Big(6^{-1/3}t^{-2/3}x/2\Big),
\end{equation}
 where
$$\kappa(x)= \frac{\mathrm{Ai}(x^2)}{x\mathrm{Ai}(x^2)+\mathrm{Ai}'(x^2)}.$$

There is an obvious similarity between Theorems \ref{Markov-LT} and \ref{mainT} and results obtained in \cite{LG1} and \cite{LGP} concerning
local times of super-Brownian motion, or equivalently of Brownian motion indexed by the Brownian tree. In particular, Theorem 1 in \cite{LGP} 
gives an equation for these local times that is analogous to \eqref{EQ1}, but with a different drift function 
(compare the preceding formula for $h(t,x)$ with the analogous formulas for $g(t,x)$ in the introduction of \cite{LGP}). Actually, the results derived in \cite{LG1} and 
\cite{LGP}  are a key ingredient of the proof of Theorem \ref{mainT}. Let us briefly explain how one can relate the process $(L^\infty_t)_{t\in\R}$
to the local times of super-Brownian motion. 

To this end, consider a super-Brownian motion $(\mathbf{X}_t)_{t\geq 0}$ with
initial value $\mathbf{X}_0=\delta_0$, and define its local times $(L_x)_{x\in\R}$ as the continuous density of the total
occupation measure $\int_0^\infty \mathbf{X}_r\,\dd r$. Also let $\dot L_x$ stand for the derivative of $L_x$ (which exists for $x\not=0$) and write $R:=\sup\{x> 0:L_x>0\}$. Then one can couple $\mathbf{X}$
and the Brownian plane $\BP$ in such a way that we have $L^\infty_t=L_{R-t}$ for all sufficiently small $t\geq 0$
(see the proof of Proposition \ref{localtimesBP} and Section \ref{Preli2} below). As a consequence of this coupling, 
the finite-dimensional marginal distributions of the process $(\lambda^{-3}L_{R-\lambda x})_{x\geq 0}$
converge when $\lambda \downarrow 0$ to the finite-dimensional marginal distributions
of $(L^\infty_x)_{x\geq 0}$ (Proposition \ref{localBP}).

An obvious idea is then to study 
the time-reversed process $(L_{R-x},\dot L_{R-x})_{x\geq 0}$ from the results of \cite{LGP}  providing a stochastic differential equation
for $(L_x,\dot L_x)_{x\geq 0}$. Unfortunately, it does not seem easy to apply the known results for the time reversal of
Markov processes to this problem.

For this reason, we follow a different route based on Proposition 18 in \cite{LGP}, which expresses 
the super-Brownian local times $(L_x)_{x\geq 0}$ in terms of a (time-changed) diffusion process $(Z_t)_{t\geq 0}$ 
solving the equation 
\begin{equation}
\label{EDS00}
\dd Z_t=4\,\dd B_t + b(Z_t)\,\dd t,
\end{equation}
where, for every $z\in\R$, 
\begin{equation}
\label{def-b}
b(z)=8\frac{p'_1(z/2)}{p_1(z/2)}-\frac{2}{3}z^2. 
\end{equation}
More precisely, we have, for every $s\geq 0$,
$$L_{\tau(s)}=L_0\,\exp\Big(\int_0^{s} Z_r\,\dd r\Big)$$
where $\tau(s)$ is the time change $\tau(s)=\inf\{x\geq 0:\int_0^x (L_y)^{-1/3}\dd y\geq s\}$. 
The process $(Z_t)_{t\geq 0}$ is recurrent and reversible with respect to its unique
invariant probability measure, which has a negative first moment (Lemma \ref{negativedrift}). One can expect to use this reversibility property to
investigate the time-reversed process $(L_{R-x})_{x\geq 0}$. Still this involves a number of
technicalities, essentially due to the fact that using the preceding formula for $L_{\tau(s)}$
requires reversing $Z$ at (random) times of the form $\sup\{s\geq 0: \int_0^s Z_r\,\dd r\geq -a\}$
(Proposition \ref{tecpro}). This line of reasoning eventually leads to Proposition \ref{repre-loc},
which shows that $L^\infty_x$ has the same distribution as $\exp(\int_0^{\tau^*_x} W^*(s)\,\dd s)$, where $(W^*_t)_{t\in\R}$ is a version (indexed by $\R$)
of the diffusion process $Z$, and $(\tau^*_x)_{x\geq 0}$ is the appropriate time change. 
One can then apply standard tools of stochastic calculus to derive equation \eqref{EQ1} in Theorem \ref{mainT}. 

It is interesting to compare the results of the present work to the paper \cite{CLG} studying hulls
in the Brownian plane. The hull of radius $r>0$ is the complement of the unique unbounded connected component
of the complement of the closed ball of radius $r$ centered at $\rho$ (informally, the hull is obtained
by filling in the bounded holes in the ball). One can define the boundary size $U_r$ of the hull of radius $r$,
and one of the main results of \cite{CLG} identifies the process $(U_r)_{r>0}$ as a time-reversed
continuous-state branching process ``starting from $+\infty$ at time $-\infty$'' and conditioned to
become extinct at time $0$. In contrast with $(U_r)_{r>0}$, the process $(L^\infty_r)_{r>0}$ is
not Markov. However, it is much smoother, since $(U_r)_{r>0}$ is only c\`adl\`ag with negative jumps (jumps correspond to times where
the hull ``swallows'' connected components of the complement of the ball). 

The paper is organized as follows. Section \ref{prelim} introduces our main objects of interest, namely
the Brownian plane, super-Brownian motion, and Brownian snake excursions. We use the construction of
the Brownian plane in \cite{CLG} 
to get a coupling between the Brownian plane and super-Brownian motion, and then to derive
Proposition \ref{localtimesBP} from known results about super-Brownian motion. Section \ref{sec:rever}
is mainly devoted to deriving certain time reversal properties and hitting distributions for the
(stationary or not stationary) solution of \eqref{EDS00}. In Section \ref{sec:repre}, we use the results of \cite{LGP}
and the preceding coupling to identify the
distribution of $(L^\infty_t,\dot L^\infty_t)_{t\geq 0}$ with that of another pair of processes that
is constructed from a solution of \eqref{EDS00} indexed by $\R$ via a time-change transformation (Proposition \ref{repre-loc}).
Finally, in Section \ref{sec:EDS}, we derive Theorems \ref{Markov-LT} and \ref{mainT} from this representation.

\section{Preliminaries}
\label{prelim}
\subsection{Brownian snake excursions and super-Brownian motion}
\label{BS-SBM}
In this section, we recall basic facts about the Brownian snake excursions and super-Brownian motion (we refer to \cite{Zurich} for more details).
We start by briefly introducing the formalism of snake trajectories.
A (one-dimensional) finite path $\w$ is a continuous mapping $\w:[0,\zeta]\la\R$, where $\zeta=\zeta_{(\w)}\in[0,\infty)$ is called the lifetime of $\w$. The space $\W$ of all finite paths is  equipped with the
distance
$$d_\W(\w,\w')=|\zeta_{(\w)}-\zeta_{(\w')}|+\sup_{t\geq 0}|\w(t\wedge
\zeta_{(\w)})-\w'(t\wedge\zeta_{(\w')})|.$$
We denote the endpoint of the path $\w$ by $\wh \w=\w(\zeta_{(\w)})$.
For every $x\in\R$, we set $\W_x=\{\w\in\W:\w(0)=x\}$. The trivial element of $\W_x$ 
with zero lifetime is identified with the point $x$ of $\R$. 

A snake trajectory with initial point $x$ is a continuous mapping $s\mapsto \omega_s$
from $\R_+$ into $\W_x$ 
which satisfies the following two properties:
\begin{enumerate}
\item[\rm(i)] We have $\omega_0=x$ and the number $\sigma(\omega):=\sup\{s\geq 0: \omega_s\not =x\}$,
called the duration of the snake trajectory $\omega$,
is finite (by convention $\sigma(\omega)=0$ if $\omega_s=x$ for every $s\geq 0$). 
\item[\rm(ii)] For every $0\leq s\leq s'$, we have
$\omega_s(t)=\omega_{s'}(t)$ for every $t\in[0,\displaystyle{\min_{s\leq r\leq s'}} \zeta_{(\omega_r)}]$ (Snake property).
\end{enumerate}

We write $\S$ for the set of all snake trajectories. If $\omega\in \S$, the occupation measure of $\omega$ is the finite measure $\mathcal{O}_\omega$ on $\R$ defined by
$$\langle \mathcal{O}_\omega,\varphi\rangle = \int_0^{\sigma(\omega)} \varphi(\wh \omega_s)\,\dd s.$$
We will write
$$\mathcal{R}(\omega)=\sup(\mathrm{supp}\,\mathcal{O}_\omega),\quad \mathcal{G}(\omega)=\inf (\mathrm{supp}\,\mathcal{O}_\omega),$$
where 
$\mathrm{supp} \,\mathcal{O}_\omega$ denotes the topological support of $\mathcal{O}_\omega$.

Let $x\in \R$. The Brownian snake excursion 
measure $\N_x$ is the $\sigma$-finite measure on $\{\omega\in\S:\omega_0=x\}$ that is characterized by the following two properties: Under $\N_x(\dd \omega)$,
\begin{enumerate}
\item[(i)] the distribution of the lifetime function $(\zeta_{\omega_s})_{s\geq 0}$ is the It\^o 
measure of positive excursions of linear Brownian motion, normalized so that $\N_x(\sup_{s\geq 0} \zeta_{\omega_s}>\ve)=\frac{1}{2\ve}$, for every $\ve>0$;
\item[(ii)] conditionally on $(\zeta_{\omega_s})_{s\geq 0}$, the function $(\wh \omega_s)_{s\geq 0}$ is
a Gaussian process with mean $x$ and covariance function 
$K(s,s')= \min_{s\wedge s'\leq r\leq s\vee s'} \zeta_{\omega_r}$.
\end{enumerate}

We record the following formula \cite[Chapter VI]{Zurich}. For every $a>x$, we have
\begin{equation}
\label{hitting-pro}
\N_x(\mathcal{R}>a)=\frac{3}{2(a-x)^2}.
\end{equation}

By results of \cite{CM} (see also \cite{BMJ} for a weaker version), the occupation measure $\mathcal{O}_\omega$ has $\N_x(\dd \omega)$ a.e.~a continuously differentiable density,
which we denote by $(\ell_y(\omega))_{y\in\R}$. Since $\mathcal{O}_\omega$ puts no mass on $(\mathcal{R}(\omega),\infty)$, both $\ell_y(\omega)$ and its derivative vanish at $y=\mathcal{R}(\omega)$.

Let us now turn to super-Brownian motion (we refer to \cite{Per} for more details about this process). We let $(\XX_t)_{t\geq 0}$ be a one-dimensional super-Brownian motion started at 
$\XX_0=\delta_0$ --- in view of the Brownian snake representation, we assume that
the branching mechanism of $\XX$ is $\Psi(\lambda)=2\lambda^2$. The ``total occupation measure'' $\mathbf{O}$ is defined by
$$\mathbf{O}=\int_0^\infty \XX_t\,\dd t$$
and the local time process $(L_x)_{x\in \R}$ is the continuous density of the random measure $\mathbf{O}$, which exists
and is continuously differentiable on $(-\infty,0)\cup (0,\infty)$, by \cite{Sug}. 
We set $R:=\sup\{x\geq 0:L_x>0\}=\sup (\mathrm{supp}\, \mathbf{O})$. By results of \cite{MP}, we have 
$L_x>0$ for every $x\in[0,R)$, a.s.

By \cite[Chapter IV]{Zurich}, we can construct the process $\XX$ from a Poisson point measure $\sum_{i\in I}\delta_{\omega_i}$
with intensity $\N_0$, in such a way that
$$\mathbf{O}=\sum_{i\in I} \mathcal{O}_{\omega_i}.$$
Write $\r(\omega_i)=\sup (\mathrm{supp}\,\mathcal{O}_{\omega_i})$. Using \eqref{hitting-pro}, we see that a.s. there is a 
unique $i_*\in I$ such that 
$$R=\sup_{i\in I} \r(\omega_i)=\r(\omega_{i_*})$$
and moreover we can find a (random) $\ve\in(0,R)$ such that
$$\sup_{i\in I\setminus\{i_*\}} \r(\omega_i)\leq R-\ve.$$
As a consequence, the restriction of $\mathbf{O}$ to $[R-\ve,R]$ coincides with 
the restriction of $\mathcal{O}_{\omega_{i_*}}$ to the same interval, and 
we have
\begin{equation}
\label{LTsnake}
L_x=\ell_x(\omega_{i_*})\,,\qquad \forall x\in[R-\ve,R].
\end{equation}

\subsection{Areas of spheres in the Brownian plane}
\label{Preli2}

The Brownian plane was introduced in \cite{plane} and further discussed in \cite{CLG} and we refer to these two papers for additional information. 
It is convenient to view the Brownian plane as a random pointed measure metric space $(\BP,D_\infty,\mathrm{Vol}, \rho)$ where $D_\infty$
is the distance on $\BP$,  $\mathrm{Vol}$ stands for the volume measure
on $\BP$ and $\rho$ is the distinguished point. The Brownian plane enjoys the following 
remarkable scale invariance property: For every $\lambda>0$,
$(\BP,\lambda D_\infty,\lambda^4\mathrm{Vol}, \rho)$ has the same distribution as $(\BP,D_\infty,\mathrm{Vol}, \rho)$.

Recall from \eqref{profi} the definition of the profile of distances $\Gamma$ in $\BP$. The scale invariance property entails 
that, for every $\lambda>0$, the scaled measure $\Gamma^{(\lambda)}$ defined by
$\Gamma^{(\lambda)}(A)=\lambda^4\Gamma(\lambda^{-1}A)$ has the same distribution as $\Gamma$. 

Let us prove Proposition \ref{localtimesBP} which was stated in Section \ref{intro}.

\proof[Proof of Proposition \ref{localtimesBP}] We rely on the construction of the Brownian plane in \cite[Section 3.2]{CLG}. This construction involves 
a nine-dimensional Bessel process $X=(X_t)_{t\geq 0}$ started at $0$, and, conditionally on $X$, two independent Poisson point measures $\n_\infty$ and $\n'_\infty$ on $\R_+\times \S$
with intensity 
$$2\,\mathbf{1}_{\{\mathcal{G}(\omega)>0\}} \dd t\,\N_{X_t}(\dd\omega),$$
in such a way that the profile of distances $\Gamma$ is the sum
of the occupation measures of all atoms of $\n_\infty$ and $\n'_\infty$
(we refer to \cite{CLG} for more details).

To make the connection with super-Brownian motion, we recall the notation introduced at the end of Section \ref{BS-SBM} 
and we  use the conditional distribution of $\omega_{i_*}$ given $\mathcal{R}(\omega_{i_*})=r$ (for any $r>0$).
This conditional distribution is $\N_0(\cdot\mid\mathcal{R}=r)$ and  is described in \cite[Theorem 2.1]{CLG}:
it involves the  nine-dimensional Bessel process $X$ 
and, conditionally on $X$, two independent Poisson point measures $\n$ and $\n'$ on $\R_+\times \S$
with intensity 
$$2\,\mathbf{1}_{[0,\tau_r]}(t)\mathbf{1}_{\{\mathcal{G}(\omega)>0\}} \dd t\,\N_{X_t}(\dd\omega).$$
where $\tau_r=\sup\{t\geq 0:X_t=r\}$. Theorem 2.1 of \cite{CLG} implies that, under the 
conditional distribution of  of $\omega_{i_*}$ given $\mathcal{R}(\omega_{i_*})=r$, the 
measure $A \mapsto \int \mathbf{1}_A(R-x)\, \mathcal{O}_{\omega_{i_*}}(\dd x)$ is equal to the sum
of the occupation measures of all atoms of $\n$ and $\n'$. 

If we compare the two representations described above
(using also Lemma 3.3 in \cite{CLG} to see that atoms of $\n_\infty$ and $\n'_\infty$ in $[1,\infty)\times \S$ do not contribute to the
values of $\Gamma$ near $0$), we easily get that  
we can couple $\omega_{i_*}$
and a Brownian plane $(\BP,D_\infty)$ in such a way that there is a
(random) $\ve'>0$ such that the restriction of the profile of distances $\Gamma$ to $[0,\ve']$ 
is equal to the restriction of the
measure $A \mapsto \int \mathbf{1}_A(R-x)\, \mathcal{O}_{\omega_{i_*}}(\dd x)$ to the same interval.
However, we know that the latter measure has a continuously differentiable density (equal to $\ell_{R-x}(\omega_{i_*})$)
on $[0,\ve']$ and that this density vanishes at $0$ together with its derivative, and is positive on $(0,\ve\wedge\ve']$ thanks to \eqref{LTsnake}. This shows that 
the profile of distances $\Gamma$ has a 
continuously differentiable density on $[0,\ve\wedge \ve']$, which vanishes at $0$ together with its derivative and is positive on $(0,\ve\wedge\ve']$. Proposition \ref{localtimesBP} now 
follows by the scale invariance of the profile of distances. \endproof

The scale invariance of the random measure $\Gamma$ readily implies a similar scale invariance property
for its density $(L^\infty_x)_{x\geq 0}$: for every $\lambda>0$ and $x\geq 0$, $L^\infty_{\lambda x}$
has the same distribution as $\lambda^3L^\infty_{x}$. We also notice that $\E[\int_0^1 L^\infty_x\dd x]=\E[\Gamma([0,1])]$
is equal to the mean volume of the unit ball in the Brownian plane, which
is finite (we can bound this volume by the volume of the hull of radius $1$, which has a  finite first moment by  \cite{CLG}). It follows that $\E[L^\infty_x]=c\,x^3$, with a constant $c\in(0,\infty)$
(one can in fact compute $c=8/21$, but we will not need this).

\begin{proposition}
\label{localBP}
The finite-dimensional marginal distributions of the process $(\lambda^{-3}L_{R-\lambda x})_{x\geq 0}$
converge when $\lambda \downarrow 0$ to the finite-dimensional marginal distributions
of $(L^\infty_x)_{x\geq 0}$.
\end{proposition}

\proof
We use the
preceding coupling obtained in the proof of Proposition \ref{localtimesBP}, under which we have
$$L^\infty_x=\ell_{R-x}(\omega_{i_*})\,,\qquad\forall x\in[0,\ve'].$$
By combining this with \eqref{LTsnake} we get
\begin{equation}
\label{couplingBP}
L^\infty_x=L_{R-x}\,,\qquad \forall x\in[0,\ve\wedge \ve'],
\end{equation}
where $\ve>0$ and $\ve'>0$ are both random. 
Then, if $0\leq x_1<x_2<\cdots<x_p$,
and $g$ is a bounded continuous function on $\R^p$, we get
$$\E[g(\lambda^{-3}L_{R-\lambda x_1},\ldots,\lambda^{-3}L_{R-\lambda x_p})]
-\E[g(\lambda^{-3}L^\infty_{\lambda x_1},\ldots,\lambda^{-3}L^\infty_{\lambda x_p})]
\build{\la}_{\lambda\to 0}^{} 0,$$
by \eqref{couplingBP} and dominated convergence. On the other hand, the scale invariance
of the Brownian plane shows that, for every $\lambda>0$,
$$\E[g(\lambda^{-3}L^\infty_{\lambda x_1},\ldots,\lambda^{-3}L^\infty_{\lambda x_p})]
=\E[g(L^\infty_{x_1},\ldots,L^\infty_{x_p})].$$
The desired result follows. \endproof

\section{The time reversal argument}
\label{sec:rever}

We consider the one-dimensional stochastic differential equation
\begin{equation}
\label{EDS}\dd Z_s= 4\,\dd B_s + b(Z_s)\,\dd s,
\end{equation}
where $B$ stands for a linear Brownian motion and  the function $b(z)$ was defined in \eqref{def-b}. From \cite[Section 6]{LGP}, we know that the solution to \eqref{EDS} is a recurrent diffusion process 
whose unique invariant probability measure is $\pi(\dd x)=\theta(x)\,\dd x$
where
$$\theta(x)=Cp_1(x/2)^2\,\exp(-\frac{x^3}{36}),$$
the function $p_1$ is given by \eqref{formupt}, 
and $C$ is the appropriate normalizing constant.
It will be convenient to introduce a collection of 
probability measures $(\P_x)_{x\in\R}$ such that, under 
$\P_x$, the process $(Z_s)_{s\geq 0}$ solves \eqref{EDS}
with initial value $Z_0=x$.

We also introduce another process $(W_s)_{s\in\R}$ indexed by $\R$,
which is distributed under the probability measure $\P$ as the stationary solution to \eqref{EDS}: for every
$t\in\R$, the process $(W_{t+s})_{s\geq 0}$ is distributed as the
solution to \eqref{EDS} with initial distribution $\pi$. It is then well known
that $(W_s)_{s\in\R}$ is time reversible, meaning
that, for every $t\in\R$, $(W_{t-s})_{s\in\R}$ has the same 
distribution as $(W_s)_{s\in\R}$ (see e.g.~\cite[Section 11]{Ken}).

\begin{lemma}
\label{negativedrift}
We have $\int_\R |x|\,\pi(\dd x) <\infty$  and
$\int_\R x\,\pi(\dd x)<0$. 
\end{lemma}

\proof
The first assertion is immediate since the function $x\mapsto p_1(x/2)$ is bounded. Let us prove the second assertion.
In the proof that follows, $c$
denotes a positive constant that may change from line to line. From formula \eqref{formupt}, we get
$$p_1(x/2)=c\,\mathcal{A}\Big(6^{-1/3}\,\frac{x}{2}\Big),$$
with $\mathcal{A}(x)=-2\,e^{2x^3/3}\,(x\mathrm{Ai}(x^2)+\mathrm{Ai}'(x^2)).$
It follows that
\begin{align*}
\int_\R x\,\pi(\dd x)&=c\int_\R x\,\mathcal{A}\Big(6^{-1/3}\,\frac{x}{2}\Big)^2\,\exp(-\frac{x^3}{36})\,\dd x\\
&=c\int_\R x\,\mathcal{A}(x)^2\,\exp(-\frac{4x^3}{3})\,\dd x\\
&=c \int_\R x\,(x\mathrm{Ai}(x^2)+\mathrm{Ai}'(x^2))^2\,\dd x
\end{align*}
Then,
\begin{align*}
&\int_\R x\,(x\mathrm{Ai}(x^2)+\mathrm{Ai}'(x^2))^2\,\dd x\\
&\qquad=\int_0^\infty x\,(x\mathrm{Ai}(x^2)+\mathrm{Ai}'(x^2))^2\,\dd x - \int_0^\infty x\,(-x\mathrm{Ai}(x^2)+\mathrm{Ai}'(x^2))^2\,\dd x\\
&\qquad=4\int_0^\infty x\,\mathrm{Ai}(x^2)\,\mathrm{Ai}'(x^2)\,\dd x
\end{align*}
which is negative since the function $\mathrm{Ai}$ is positive and monotone decreasing on $[0,\infty)$. \endproof

From Lemma \ref{negativedrift} and the ergodic theorem we have $\P$ a.s.,
$$\lim_{t\to\infty} \int_0^t W_u\,\dd u=-\infty\;,\quad \lim_{t\to-\infty} \int_t^0 W_u\,\dd u=-\infty.$$

We set
\begin{align*}
D&:=\Big\{s\in\R: \forall t\in(-\infty,s), \int_t^s W_u\,\dd u<0\Big\},\\
D^*&:=\Big\{s\in\R: \forall t\in(s,\infty), \int_s^t W_u\,\dd u<0\Big\}.
\end{align*}
Clearly, we have 
\begin{equation}
\label{negativity}
W_s\leq 0,\hbox{ for every }s\in D\cup D^*.
\end{equation}
For every $x\in\R$, we also set
\begin{equation}
\label{def-gamma}
\gamma(x):=\P_x\Big(\forall t>0, \int_0^t Z_u\,\dd u<0\Big).
\end{equation}
It is then easy to verify that 
\begin{equation}
\label{gamma-positive}
\gamma(x)>0 \hbox{ if and only if }x<0.
\end{equation}

In the remaining part of this section, {\bf we fix} $a>0$
(many of the random quantities that we will introduce will depend on $a$, 
although this will not be apparent in the notation). For every $t\in\R$, we set
$$T_t:=\sup\Big\{s\in[t,\infty): \int_t^s W_u\,\dd u\geq -a\Big\}.$$
By construction, $\P$ a.s. for every $t\in\R$, we have $T_t<\infty$, $\int_t^{T_t} W_u\,\dd u=-a$
and, for every $r>T_t$, $\int_t^{r} W_u\,\dd u<-a$, so that
$\int_{T_t}^rW_u\,\dd u<0$, and thus $T_t\in D^*$. In particular
we must have $W_{T_t}\leq 0$ by \eqref{negativity}. 

Similarly, we set, for every $t\in\R$,
$$T^*_t=\inf\Big\{s\in(-\infty,t]: \int_s^t W_u\,\dd u\geq -a\Big\}.$$
Then $T^*_t\in D$ and we have again $W_{T^*_t}\leq 0$. 
The mapping $D\ni t\mapsto T_t\in D^*$ is a bijection from $D$ onto $D^*$,
whose inverse is $D^*\ni s\mapsto T^*_s\in D$. 

We then observe that the set $\wt D:=\{t\in D: W_t<0\hbox{ and } W_{T_t}<0\}$ is open. Indeed,
suppose  that $s\in \wt D$. Then, for $\ve>0$ small enough
(such that $W_r<0$ for $r\in[s,s+\ve]$), it is immediate that $\int_t^{s+\ve} W_u\,\dd u<0$
for every $t<s+\ve$, so that $s+\ve\in D$ and $W_{s+\ve}<0$. Let us then verify that $T_{s+\ve}\to T_s$ as $\ve\downarrow 0$.
Since $W_s<0$, the mapping $r\mapsto\int_s^r W_u\,\dd u$ is monotone decreasing
on a small interval $[s,s+\delta]$, $\delta>0$, and, by construction, this implies that $\ve\mapsto T_{s+\ve}$
is monotone decreasing on $[0,\delta]$. Therefore, $T_{s+\ve}$ has a limit 
when $\ve\downarrow 0$, $\ve>0$, and we denote this limit by $T_{s+}\geq T_s$. However, since $\int_{s+\ve}^{T_{s+\ve}} W_u\,\dd u=-a$, we immediately
get that $\int_s^{T_{s+}} W_u\,\dd u=-a$, so that we have also $T_{s+}\leq T_s$ and we have proved our
claim $T_{s+}= T_s$. Since $W_{T_s}<0$, it follows that
we have $W_{T_{s+\ve}}<0$ for $\ve>0$ small enough, and therefore, as we already noticed that 
$s+\ve\in D$ and $W_{s+\ve}<0$, we get that
$s+\ve \in \wt D$ if $\ve$ is small enough. A similar argument, left
to the reader, shows that $s-\ve\in \wt D$ if $\ve$ is small enough, and
we have proved that $\wt D$ is open. 

We can also set $\wt D^*:=\{s\in D^*: W_s<0\hbox{ and } W_{T^*_s}<0\}$. A symmetric argument shows that
$\wt D^*$ is open and the mapping $t\mapsto T_t$ is a bijection from $\wt D$ onto $\wt D^*$. 

Suppose that $t\in \wt D$. Then 
$$\int_t^{t+\ve} W_u\,\dd u\build{\sim}_{\ve \to 0}^{} \ve\,W_t$$
and
$$\int_{T_t}^{T_{t+\ve}} W_u\,\dd u \build{\sim}_{\ve \to 0}^{} (T_{t+\ve}-T_\ve)W_{T_t}$$
and since $\int_t^{T_t} W_u\,\dd u=-a=\int_{t+\ve}^{T_{t+\ve}} W_u\,\dd u$, it follows that
$$(T_{t+\ve}-T_t)\,W_{T_t} \build{\sim}_{\ve \to 0}^{}  \ve \,W_t.$$
Hence the mapping $r\mapsto T_r$ is differentiable at $t$, and its 
derivative is $W_t/W_{T_t}$.

Let $g$ be a measurable function from $\R_+$ into $\R_+$, let $\Phi$ be a nonnegative
measurable function on the space $\mathcal{W}$ of all finite paths, and let
$F$ be a nonnegative measurable function on the space $C(\R_+,\R)$
of all continuous functions from $\R_+$ into $\R_+$. The change of variables $s=T_t$
($t=T^*_s$) shows that
\begin{align}
\label{changevar}
&\qquad\E\Big[\int_{-\infty}^{\infty} \dd t\,g(t)\,|W_t|\,\mathbf{1}_{\wt D}(t)\,\Phi (W_{T_t-u}:0\leq u\leq T_t-t)\,F(W_{T_t+r}:r\geq 0)\Big]\nonumber\\
&\qquad=\E\Big[\int_{-\infty}^{\infty} \dd s\,g(T^*_s)\,|W_s|\,\mathbf{1}_{\wt D^*}(s)\,\Phi (W_{s-u}:0\leq u\leq s-T^*_s)\,F(W_{s+r}:r\geq 0)\Big].
\end{align}
We can in fact replace $\mathbf{1}_{\wt D}(t)$ by $\mathbf{1}_{D}(t)$ in the LHS of \eqref{changevar}, and similarly $\mathbf{1}_{\wt D^*}(s)$
by $\mathbf{1}_{D^*}(s)$ in the RHS of \eqref{changevar}. Recalling \eqref{negativity}, the point is to observe that we have
 a.s.
$$\int_{-\infty}^{\infty} \mathbf{1}_{\{W_t=0\}}\,\dd t=0\hbox{ and } \int_{-\infty}^{\infty} \mathbf{1}_{\{W_{T_t}=0\}}\,\dd t=0,$$
where the second equality can be derived by an absolute continuity argument and properties of
linear Brownian motion showing that, for every fixed $t$, there are a.s. no $s>t$ such that $\int_t^s W_r\,\dd r=-a$ and $W_s=0$ (see e.g.~Lachal \cite{Lac}). 

Suppose that $\Phi$ has the additional property that $\Phi(\w)=0$ if $\zeta_{(\w)}>M$, for some $M>0$, and take
$g=\mathbf{1}_{[-A,A]}$, where $A>0$. By stationarity, the LHS of \eqref{changevar} (with $\mathbf{1}_{\wt D}(t)$ replaced by $\mathbf{1}_{D}(t)$)
is equal to
$$2A\,\E\Big[|W_0|\,\mathbf{1}_D(0)\,\Phi (W_{T_0-u}:0\leq u\leq T_0)\,F(W_{T_0+r}:r\geq 0)\Big].$$
On the other hand, our assumption on $\Phi$ implies that, if $\Phi (W_{s-u}:0\leq u\leq s-T^*_s)\not=0$, we have $s-M\leq T^*_s\leq s$ and thus
$\mathbf{1}_{[-A,A]}(T^*_s)$ is bounded above by $\mathbf{1}_{[-A,A+M]}(s)$ and bounded below by $\mathbf{1}_{[-A+M,A]}(s)$. Hence, the 
RHS of \eqref{changevar} is bounded above by
$$(2A+M)\,\E\Big[|W_0|\,\mathbf{1}_{ D^*}(0)\,\Phi (W_{-u}:0\leq u\leq -T^*_0)\,F(W_{r}:r\geq 0)\Big]$$
and bounded below by the same quantity with $(2A+M)$ replaced by $(2A-M)$. If we multiply the resulting bounds by $1/2A$ and let $A\to\infty$,
we arrive at the formula
\begin{align}
\label{keyform}
&\qquad\E\Big[|W_0|\,\mathbf{1}_D(0)\,\Phi (W_{T_0-u}:0\leq u\leq T_0)\,F(W_{T_0+r}:r\geq 0)\Big]\nonumber\\
&\qquad\qquad=\E\Big[|W_0|\,\mathbf{1}_{D^*}(0)\,\Phi (W_{-u}:0\leq u\leq -T^*_0)\,F(W_{r}:r\geq 0)\Big].
\end{align}
Clearly, this formula remains valid without the extra assumption we made on $\Phi$ (replace $\Phi(\w)$
by $\Phi(\w)\,\mathbf{1}_{\{\zeta_{(\w)}\leq n\}}$ and let $n\to\infty$).
Taking $F=1$ and $\Phi(\w)=\varphi(\w(0))$ for a nonnegative measurable function $\varphi:\R\la \R_+$, we get
$$\E\Big[|W_0|\,\mathbf{1}_D(0)\,\varphi(W_{T_0})\Big]
=\E\Big[|W_0|\,\mathbf{1}_{D^*}(0)\,\varphi(W_0)\Big].$$
Note that the Markov property and the invariance of $W$ under time-reversal imply that
$\P(0\in D\,|\, (W_r)_{r\geq 0})=\P(0\in D^*\,|\, (W_r)_{r\leq 0})=\gamma(W_0)$, where $\gamma$ was defined in \eqref{def-gamma}.  It follows that
$$\E\Big[|W_0|\,\gamma(W_0)\,\varphi(W_{T_0})\Big]
=\E\Big[|W_0|\,\gamma(W_0)\,\varphi(W_0)\Big].$$
Hence, under the (suitably normalized) initial distribution $|y|\gamma(y)\,\pi(\dd y)$, $W_{T_0}$ has the same distribution as $W_0$. 

We set $\mu(\dd y)=c_0\,|y|\gamma(y)\,\pi(\dd y)$, where the constant $c_0>0$ is chosen so that $\mu$ is a probability measure.
Notice that $\mu$ is supported on the negative half-line, by \eqref{gamma-positive}. 
We introduce the probability measure $\P_\mu$ under which the process $(Z_s)_{s\geq 0}$ is 
a solution to \eqref{EDS} with initial distribution $\mu$. We also set 
$$\xi=\sup\Big\{s\in[0,\infty): \int_0^s Z_u\,\dd u\geq -a\Big\},$$ and
write $\Xi$
for the event
$$\Xi=\Big\{\int_0^t Z_s\,\dd s<0,\;\forall t>0\Big\}.$$

\begin{proposition}
\label{tecpro}
The following properties hold under $\P_\mu$.

{\rm1.} The two processes $(Z_{\xi-u})_{0\leq u\leq \xi}$ and $(Z_{u})_{0\leq u\leq \xi}$ have the same distribution.
In particular $Z_{\xi}$ is distributed according to $\mu$.

{\rm 2.} The processes $(Z_{\xi-u})_{0\leq u\leq \xi}$ and $(Z_{\xi+r})_{r\geq 0}$ are conditionally independent
given $Z_{\xi}$, and, for every $x<0$, the conditional distribution of $(Z_{\xi+r})_{r\geq 0}$ knowing $Z_{\xi}=x$
is the distribution of $(Z_r)_{r\geq 0}$ under $\P_x(\cdot\mid \Xi)$.

{\rm 3.} For $b>0$, let $H_b:=\inf\{t\geq 0: \int_0^t Z_r\,\dd r=-b\}$. Then $Z_{H_b}$ is distributed according to $\mu$.
\end{proposition}

\proof 1. By taking $F=1$ in \eqref{keyform}, and then conditioning on $(W_r)_{r\geq 0}$ in the LHS and on
$(W_r)_{r\leq 0}$ in the RHS, we arrive at
\begin{equation}
\label{keyform2}
\E\Big[|W_0|\,\gamma(W_0)\,\Phi (W_{T_0-u}:0\leq u\leq T_0)\Big]=\E\Big[|W_0|\,\gamma(W_0)\,\Phi (W_{-u}:0\leq u\leq -T^*_0)\Big].
\end{equation}
However, since $(W_s)_{s\in\R}$ and $(W_{-s})_{s\in\R}$ have the same distribution, one easily gets that
$$\E\Big[|W_0|\,\gamma(W_0)\,\Phi (W_{-u}:0\leq u\leq -T^*_0)\Big]=\E\Big[|W_0|\,\gamma(W_0)\,\Phi (W_{u}:0\leq u\leq T_0)\Big],$$
and thus
$$\E\Big[|W_0|\,\gamma(W_0)\,\Phi (W_{T_0-u}:0\leq u\leq T_0)\Big]=\E\Big[|W_0|\,\gamma(W_0)\,\Phi (W_{u}:0\leq u\leq T_0)\Big].$$
The desired result follows since, by construction, the distribution of $((Z_s)_{s\geq 0},\xi)$ under $\P_\mu$
coincides with the distribution of $((W_s)_{s\geq 0},T_0)$ under the probability measure which has density $c_0|W_0|\gamma(W_0)$
with respect to $\P$. 

2. Similarly, we have
\begin{align*}
&\E_\mu\Big[\Phi (Z_{\xi-u}:0\leq u\leq \xi)\,F(Z_{\xi+r}:r\geq 0)\Big]\\
&\qquad=c_0\,\E\Big[|W_0|\,\gamma(W_0)\,\Phi (W_{T_0-u}:0\leq u\leq T_0)\,F(W_{T_0+r}:r\geq 0)\Big]\\
&\qquad=c_0\,\E\Big[|W_0|\,\mathbf{1}_D(0)\,\Phi (W_{T_0-u}:0\leq u\leq T_0)\,F(W_{T_0+r}:r\geq 0)\Big]\\
&\qquad=c_0\,\E\Big[|W_0|\,\mathbf{1}_{D^*}(0)\,\Phi (W_{-u}:0\leq u\leq -T^*_0)\,F(W_{r}:r\geq 0)\Big]\\
&\qquad=c_0\,\E\Big[|W_0|\,\gamma(W_0)\,\Phi (W_{-u}:0\leq u\leq -T^*_0)\,\E_{W_0}[F(Z_{r}:r\geq 0)\mid {\Xi}\,]\Big]\\
&\qquad=c_0\,\E\Big[|W_0|\,\gamma(W_0)\,\Phi (W_{T_0-u}:0\leq u\leq T_0)\,\E_{W_{T_0}}[F(Z_{r}:r\geq 0)\mid {\Xi}\,]\Big]\\
&\qquad=\E_\mu\Big[\Phi (Z_{\xi-u}:0\leq u\leq \xi)\,\E_{Z_{\xi}}[F(Z_{r}:r\geq 0)\mid {\Xi}\,]\Big].
\end{align*}
We used \eqref{keyform} in the third equality, then the Markov property of $W$ at time $0$, and \eqref{keyform2} in the fifth equality.

3. We may assume that $b<a$. We then observe that
$$\xi-H_b=\sup\Big\{s\in [0,\xi]: \int_0^s Z_{\xi-r}\,\dd r\geq b-a\Big\}.$$
Writing $Z'(u)=Z_{\xi-u}$ for $0\leq u\leq \xi$, it follows that
$$Z_{H_b}=Z'(\xi-H_b)=Z'\Big(\sup\Big\{s\in [0,\xi]: \int_0^s Z'_r\,\dd r\geq b-a\Big\}\Big),$$
which is distributed according to $\mu$ by part 1 of the proposition. 
\endproof

Suppose now that $Z$ is a solution to \eqref{EDS} with an arbitrary initial distribution, and keep the notation
$H_b=\inf\{t\geq 0:\int_0^t Z_r\,\dd r=-b\}$. In view of future applications, we would like to say that the distribution of $Z_{H_b}$ is close to $\mu$
in total variation when $b$ is large. We content ourselves with a slightly weaker result.

\begin{proposition}
\label{couplingTV}
Let $Z$ be a solution to \eqref{EDS} and let $\ve>0$.
The following holds for every large enough $c>0$. If $\kappa_c$ is a random variable independent of $Z$
and uniformly distributed over $[c,2c]$, the total variation distance between $\mu$ and the
distribution of $Z_{H_{\kappa_c}}$ is less than $\ve$. 
\end{proposition}

\proof We use a standard coupling argument (see \cite[Lemma 23.17]{Kal}). We may consider a process $\wh Z$,
which solves \eqref{EDS} with initial distribution $\mu$, such that $Z$ and $\wh Z$ are coupled 
in the following way: The random variable $T:=\inf\{t\geq 0:Z_t=\wh Z_t\}$
is a.s. finite, and $Z_t=\wh Z_t$ for every $t\geq T$. We then choose a constant $M>0$
such that the probability of the event
$$A_M:=\Big\{T\leq M, \int_0^T|Z_s|\,\dd s\leq M,  \int_0^T|\wh Z_s|\,\dd s\leq M\Big\}$$
is at least $1-\ve/4$. 

Set $\wh H_c=\inf\{t\geq 0:\int_0^t \wh Z_r\,\dd r=-c\}$, and suppose that $c$ is large enough so that the probability of the event 
$$B_M:=\{H_c\geq2M,\wh H_c\geq 2M\}$$
is at least $1-\ve/4$. Then, on the event $A_M\cap B_M$, for every $u>4M$,
we have $\wh H_u>T$ and thus
\begin{align*}
\wh H_u=\inf\Big\{t\geq T:\int_0^t \wh Z_s\,\dd s=-u\Big\}
&=\inf\Big\{t\geq T:\int_0^t Z_s\,\dd s=-u+\int_0^T Z_s\,\dd s-\int_0^T\wh Z_s\,\dd s\Big\}\\
&=H_{u-X},
\end{align*}
where $X:=\int_0^T Z_s\,\dd s-\int_0^T\wh Z_s\,\dd s$, and we note that the condition $u>4M$ ensures that $H_{u-X}>T$.

Let $\varphi:\R\la [0,1]$ be a Borel function. Still on the event $A_M\cap B_M$, we get, for $c>4M$,
$$\int_c^{2c} \Big(\varphi(Z_{H_u})-\varphi(\wh Z_{\wh H_u})\Big)\,\dd u
=\int_c^{2c} \Big(\varphi(Z_{H_u})-\varphi(Z_{H_{u-X}})\Big)\,\dd u\leq 2|X|\leq 4M.$$
Finally (as we may assume that $\kappa_c$ is also independent of $\wh Z$),
we get
$$|\E[\mathbf{1}_{A_M\cap B_M}(\varphi(Z_{H_{\kappa_c}})-\varphi(\wh Z_{\wh H_{\kappa_c}}))]|
\leq \frac{4M}{c}<\frac{\ve}{2}$$
when $c$ is large, independently of the choice of $\varphi$. Since $\P(A_M\cap B_M)\geq 1-\ve/2$, it follows that
$$|\E[(\varphi(Z_{H_{\kappa_c}})]-\E[\varphi(\wh Z_{\wh H_{\kappa_c}})]|\leq \ve,$$
which completes the proof since the distribution of $\wh Z_{\wh H_{\kappa_c}}$ is $\mu$ by the preceding proposition. \endproof

\rem We can slightly extend Proposition \ref{couplingTV} as follows. Suppose that $Z$ is Markov
with respect to a filtration $(\g_t)_{t\geq 0}$, and let $Y$ be a $\g_0$-measurable random variable. Then, the distribution of
$Z_{H_{Y+\kappa_c}}$ will again be close to $\mu$ in total variation when $c$ is large. This follows by a minor
modification of the preceding proof. 

\section{A representation for the process of local times}
\label{sec:repre}

Recall that $(L_x)_{x\in\R}$ are the local times of the super-Brownian motion $\mathbf{X}$ started at $\delta_0$, and
$R=\sup\{x\geq 0: L_x>0\}$. We have $L_x>0$ for every $x\in[0,R)$. We write $\dot L_x$ for the derivative of
$L_x$ (when $x=0$, we define $\dot L_0$ as the right derivative at $0$, see \cite{Sug}).

For every $t\geq 0$,
set
$$\tau(t):= \inf\Big\{x\geq 0: \int_0^x (L_y)^{-1/3}\,\dd y\geq t\Big\},$$
which makes sense and belongs to $[0,R)$ because
$$\int_0^R (L_y)^{-1/3}\,\dd y=\infty,$$
by \cite[Proposition 18]{LGP}.  We set $\wt L(t):=L_{\tau(t)}$ for every $t\geq 0$.
Then, according to \cite[Proposition 18]{LGP}, we have
\begin{equation}
\label{defLtilde}
\wt L(t)=\wt L(0)\,\exp\Big(\int_0^t \wt Z(s)\,\dd s\Big),
\end{equation}
where the process $(\wt Z(s))_{s\geq 0}$ solves the equation \eqref{EDS} 
with initial value $\wt Z(0)=\dot L_0/(L_0)^{2/3}$. We also note that, for $t\in [0,\infty)$,
$$\tau(t)=\int_0^t (\wt L(s))^{1/3}\,\dd s,$$
and in particular, by letting $t\to\infty$,
$$R=\int_0^\infty (\wt L(s))^{1/3}\,\dd s.$$
The inverse of $\tau(\cdot)$ is given for $x\in[0,R)$ by
\begin{equation}
\label{inverse-tau}
\tau^{-1}(x)=\inf\Big\{t:\int_0^t (\wt L(s))^{1/3}\dd s\geq x\Big\}.
\end{equation}
Then, for $x\in(0,R]$, we have using \eqref{inverse-tau}
\begin{equation}
\label{repre1}
L_{R-x}=\wt L\Big(\inf\Big\{t:\int_0^t (\wt L_s)^{1/3}\dd s\geq R-x\Big\}\Big)
=\wt L\Big(\sup\Big\{t:\int_t^\infty (\wt L_s)^{1/3}\dd s\geq x\Big\}\Big). 
\end{equation}

The next lemma introduces a process $W^*$ which plays a key role in what follows. 
\begin{lemma}
\label{constructW}
We can construct 
 a process $(W^*(t))_{t\in\R}$
with continuous sample paths, whose law is characterized by the following properties. For every $b\geq 0$, set
\begin{equation}
\label{defSb}
S_b:=\inf\Big\{t\leq 0: \int_t^0 W^*(s)\,\dd s\geq -b\Big\}.
\end{equation}
Then, $S_0=0$ a.s. Moreover, for every $b\geq0$, $S_b>-\infty$ a.s. and the process $(W^*(S_b+t))_{t\geq 0}$ is distributed as $(Z_t)_{t\geq 0}$ under $\P_\mu$,
\end{lemma}

\proof Let us first fix $b\geq 0$. Let $(Z'_t)_{t\geq0}$ be distributed as the solution of \eqref{EDS}
with initial distribution $\mu$. Set
$$\Sigma_b:=-\inf\Big\{t\geq 0:\int_0^t Z'_s\,\dd s=-b\Big\},$$
and define a process $(V_b(t),t\in[\Sigma_b,\infty))$ by setting
$V_b(t)=Z'_{t-\Sigma_b}$ for every $t\geq \Sigma_b$. Notice that we can also write
\begin{equation}
\label{cons1}
\Sigma_b=-\inf\Big\{t\geq 0:\int_0^t V_b(\Sigma_b+s)\,\dd s=-b\Big\},
\end{equation}
and that $\int_{\Sigma_b}^0 V_b(s)\,\dd s=-b$. 
Then, suppose that $0\leq b'<b$, and
set
$$\wt\Sigma_{b'}=\inf\Big\{s\geq \Sigma_b:\int_{\Sigma_b}^s V_b(r)\,\dd r =b'-b\Big\}.$$
so that $\Sigma_b<\wt \Sigma_{b'}\leq 0$. Observe that, for every $t\geq 0$,
$$V_b(\wt\Sigma_{b'}+t)=Z'_{\wt\Sigma_{b'}-\Sigma_b+t},$$
with 
$$\wt\Sigma_{b'}-\Sigma_b=\inf\Big\{s\geq 0:\int_0^s V_b(\Sigma_b+r)\,\dd r=b'-b\Big\}= \inf\Big\{s\geq 0:\int_0^s Z'_r\,\dd r=b'-b\Big\}.$$
It then follows from the Markov property of $Z'$ and point 3. of Proposition \ref{tecpro} that the 
process $(V_b(\wt\Sigma_{b'}+t))_{t\geq 0}$ is distributed as the solution of \eqref{EDS}
with initial distribution $\mu$. Furthermore, elementary considerations show that we have also
\begin{equation}
\label{cons2}
\wt\Sigma_{b'}=-\inf\Big\{t\geq 0:\int_0^t V_b(\wt\Sigma_{b'}+s)\,\dd s=-b'\Big\}.
\end{equation}

The preceding observations imply that $(V_b(\wt\Sigma_{b'}+t))_{t\geq 0}$ and $(V_{b'}(\Sigma_{b'}+t))_{t\geq 0}$ have
the same distribution, and, more precisely (using \eqref{cons1} and \eqref{cons2}), the pairs 
$((V_b(\wt\Sigma_{b'}+t))_{t\geq 0},\wt\Sigma_{b'})$ and $((V_{b'}(\Sigma_{b'}+t))_{t\geq 0},\Sigma_{b'})$ have the same
distribution. 

By taking a projective limit (or using Kolmogorov's extension theorem), we can construct
simultaneouly all processes $(V_b(t),t\in[\Sigma_b,\infty))$ for $b\geq 0$ in such a way that,
if $0\leq b'<b$, we have both
\begin{equation}
\label{cons3}
\Sigma_{b'}=\inf\Big\{s\geq \Sigma_b:\int_{\Sigma_b}^s V_b(r)\,\dd r =b'-b\Big\},
\end{equation}
and 
$$V_{b'}(t)=V_b(t)\hbox{ for every }t\geq \Sigma_{b'}.$$
Furthermore, it is easy to verify that we have then $\Sigma_b\la -\infty$ as $b\to\infty$ (writing $-\Sigma_n=\Sigma_{n-1}-\Sigma_n+
\Sigma_{n-2}-\Sigma_{n-1}+\cdots+\Sigma_1-\Sigma_0$ shows that the limit of $-\Sigma_b$ at $+\infty$
is the sum of an infinite sequence of identically distributed positive random variables). This makes it possible 
to define, for every $t\in\R$, $W^*(t)=V_b(t)$, independently of the choice of $b$
such that $\Sigma_b\leq t$. Moreover, using \eqref{cons1} and \eqref{cons3}, one verifies that the random variable 
$S_b$ defined in the statement coincides with $\Sigma_b$, so that the various
properties of the process $W^*$ follow easily. 

Finally, to see that these properties characterize the distribution of $W^*$, suppose that
$\check W$ is a process satisfying the same properties, and let $\check S_b$ be defined by
replacing $W^*$ by $\check W$ in the definition of $S_b$. The property $\check S_0=0$
ensures that $\int_t^0 \check W(s)\,\dd s<0$ for every $t<0$. Using this property, it is then straightforward to
verify that
$$\check S_b=-\inf\Big\{t\geq 0:\int_0^t \check W(\check S_b+s)\,\dd s=-b\Big\}.$$
Comparing with \eqref{cons1} (with $V_b$ and $\Sigma_b$ replaced by
$W^*$ and $S_b$ respectively), we get that the pair $((\check W(\check S_b+t))_{t\geq 0}, \check S_b)$
has the same law as $((W^*(S_b+t))_{t\geq 0}, S_b)$, and since this holds for any
$b\geq 0$, it follows that $\check W$ and $W^*$ have the same distribution. \endproof

In what follows, $W^*$ denotes the process introduced in Lemma \ref{constructW} and we keep
the notation $S_b$ introduced in \eqref{defSb}. We set, for every $t\in\R$,
$$\Lambda^*(t)=\exp\Big(\int_0^t W^*(s)\,\dd s\Big).$$
Note that $t^{-1}\int_0^t W^*_s\dd s$ converges a.s.~to $\int x\,\pi(dx)<0$ as $t\to+\infty$, and in particular
$\Lambda^*(t)$ tends to $0$ exponentially fast as $t\to+\infty$, a.s. On the other hand, 
the fact that $S_b>-\infty$ for every $b\geq 0$ ensures that $\int_0^t W^*(s)\,\dd s\la +\infty$, and therefore $\Lambda^*(t)\la +\infty$,
as $t\to-\infty$.

We set, for every $x>0$,
$$\tau^*_x=\sup\Big\{z\in \R: \int_z^\infty \Lambda^*(s)^{1/3}\,\dd s\geq x\Big\}.$$
Then $x\mapsto \tau^*_x$
is a monotone decreasing continuous function mapping $(0,\infty)$ onto $\R$. 

\begin{proposition}
\label{repre-loc}
The process $(L^\infty_x)_{x>0}$ has the same distribution as the process
$(\Lambda^*(\tau^*_x))_{x>0}$. More precisely, $(L^\infty_x,\dot L^\infty_x)_{x>0}$
has the same distribution as $(\Lambda^*(\tau^*_x), -\Lambda^*(\tau^*_x)^{2/3}W^*(\tau^*_x))_{x>0}$. 
\end{proposition}

\proof For every $\lambda>0$, let $\theta_\lambda$ be a positive random variable independent of the super-Brownian motion $\mathbf{X}$ 
and such that $-\log\,\theta_\lambda$ is uniformly distributed over $[\lambda,2\lambda]$. 
From Proposition \ref{localBP}, we get that the finite-dimensional marginal distributions of the process 
$(\theta_\lambda^{-3}L_{R-\theta_\lambda x})_{x\geq 0}$ converge to those of $(L^\infty_x)_{x\geq 0}$
when $\lambda\to\infty$. We will use \eqref{repre1} to study $(\theta_\lambda^{-3}L_{R-\theta_\lambda x})_{x\geq 0}$  when $\lambda\to\infty$. 
We fix $\ve >0$ and $M>0$. We first choose $K>1$ large enough so that the 
bound
\begin{equation}
\label{repre0}
K^{1/3}\int_0^\infty \Lambda^*(s)^{1/3}\,\dd s >M
\end{equation}
holds with probability at least $1-\ve$. 

By \eqref{repre1},
we have for every $x>0$ and $\lambda>0$ such that $\theta_\lambda x<R$,
\begin{equation}
\label{repre00}
\theta_\lambda^{-3}L_{R-\theta_\lambda x}=\theta_\lambda^{-3}\,\wt L\Big(\sup\Big\{t:\int_t^\infty (\theta_\lambda^{-3}\wt L(s))^{1/3}\dd s\geq x\Big\}\Big).
\end{equation}
Define
$$c_{\lambda,K}:=\inf\{t\geq 0: \theta_\lambda^{-3}\wt L(t)=K\}$$
where $\inf\varnothing=\infty$. We have $c_{\lambda,K}<\infty$ on the event $\{\wt L(0)\geq K \theta_\lambda^3\}$, and (since $\wt L(0)=L_0>0$ a.s.) we may assume that $\lambda>0$ is large enough so that
the probability of the event $\{\wt L(0)\geq \theta_\lambda^3K\}$ is at least $1-\ve$. From
now on we argue on this event. By \eqref{defLtilde}, we have for every $r\geq 0$,
\begin{equation}
\label{repre11}
\wt L(c_{\lambda,K}+r)= \wt L(c_{\lambda,K})\,\exp\Big(\int_{c_\lambda}^{c_\lambda+r} \wt Z(s)\,\dd s\Big)
= K\theta_\lambda^3\,\exp\Big(\int_0^r \wt Z(c_{\lambda,K}+s)\,\dd s\Big).
\end{equation}
Then (conditionally on the event $\{\wt L(0)\geq K\theta_\lambda^3\}$), we know that $(\wt Z(c_{\lambda,K}+s))_{s\geq 0}$ solves
\eqref{EDS} with initial value
\begin{align*}
\wt Z(c_{\lambda,K})&=\wt Z\Big(\inf\Big\{t\geq 0:\theta_\lambda^{-3}\wt L(0) \exp\Big(\int_0^t\wt Z(s)\dd s\Big)=K\Big\}\Big)\\
&=\wt Z\Big(\inf\Big\{t\geq 0:\int_0^t\wt Z(s)\dd s=\log(K/\wt L(0))+3\log\theta_\lambda\Big\}\Big).
\end{align*}
Notice that, by construction, $-3\log\theta_\lambda$ is uniformly distributed over $[3\lambda,6\lambda]$, and so we 
know from Proposition \ref{couplingTV} and the subsequent remark that the distribution of $\wt Z(c_{\lambda,K})$ is close to $\mu$
in total variation when $\lambda$ is large. 
It follows that, for $\lambda$ large, we can couple $(\wt Z(c_{\lambda,K}+s))_{s\geq 0}$ with $(W^*(s))_{s\geq 0}$
in such a way that we have 
$$\wt Z(c_{\lambda,K}+s)=W^*(s)$$
for every $s\geq 0$, except possibly on an event of probability at most $\ve$. Under this coupling, and using \eqref{repre11}, we get that, except on an event of probability at most $2\ve$, we
have both $c_{\lambda,K}<\infty$ and, for every $r\geq 0$,
\begin{equation}
\label{repre2}
\theta_\lambda^{-3}\wt L(c_{\lambda,K}+r)=K\exp\Big(\int_0^r W^*(s)\,\dd s\Big)=K\,\Lambda^*(r).
\end{equation}
Discarding another event of probability at most $\ve$, we know that \eqref{repre0} holds, 
so that 
\begin{equation}
\label{tec23}
\int_0^\infty (K\Lambda^*(s))^{1/3}\,\dd s >M
\end{equation}
and then,  by \eqref{repre2},
$$\int_{c_{\lambda,K}}^\infty (\theta_\lambda^{-3}\wt L(s))^{1/3}\dd s >M,$$
so that, for every $x\in(0,M]$,
\begin{align*}
\sup\Big\{t\geq 0:\int_t^\infty (\theta_\lambda^{-3}\wt L(s))^{1/3}\dd s\geq x\Big\}&=c_{\lambda,K}+\sup\Big\{z\geq 0:\int_{c_{\lambda,K}+z}^\infty (\theta_\lambda^{-3}\wt L(s))^{1/3}\dd s\geq x\Big\}\\
&=c_{\lambda,K}+\sup\Big\{z\geq 0:\int_z^\infty (K\Lambda^*(s))^{1/3}\,\dd s \geq x\Big\},
\end{align*}
using \eqref{repre2} again. From \eqref{repre00}, the last display and  \eqref{repre2}, we get, for $x\in(0,M]$,
\begin{align*}
\theta_\lambda^{-3}L_{R-\theta_\lambda x}&=\theta_\lambda^{-3}\,\wt L\Big(c_{\lambda,K}+\sup\Big\{z\geq 0:\int_z^\infty (K\Lambda^*(s))^{1/3}\,\dd s \geq x\Big\}\Big)\\
&=K\,\Lambda^*\Big(\sup\Big\{z\geq 0:\int_z^\infty (K\Lambda^*_s)^{1/3}\,\dd s \geq x\Big\}\Big).
\end{align*}
and we can replace $z\geq 0$ by $z\in\R$ in the last line because of \eqref{tec23}.

Set $b=\log K$ and recall the definition of $S_b$ in \eqref{defSb}.  Observe that $\Lambda^*(S_b+t)=K\,\exp(\int_{S_b}^{S_b+t} W^*(s)\,\dd s)$, from which it follows that the process $(K\Lambda^*(t))_{t\in\R}$ 
has the same distribution as $(\Lambda^*(S_b+t))_{t\in \R}$. It easily follows that the process
$$\Big(K\,\Lambda^*\Big(\sup\Big\{z\in\R:\int_z^\infty (K\Lambda^*_s)^{1/3}\,\dd s \geq x\Big\}\Big)\Big)_{x\in (0,M]}$$
has the same distribution as the process $(\Lambda^*(\tau^*_x))_{x\in(0,M]}$. 

Summarizing, for every large enough $\lambda$, except on a set of probability at most $3\ve$, $(\theta_\lambda^{-3}L_{R-\theta_\lambda x})_{x\in(0,M]}$ coincides
with a process which has the same finite-dimensional marginals as the process $(\Lambda^*(\tau^*_x))_{x\in(0,M]}$.
Since the finite-dimensional marginals of $(\theta_\lambda^{-3}L_{R-\theta_\lambda x})_{x\geq 0}$ converge to those of $(L^\infty_x)_{x\geq 0}$
when $\lambda\to\infty$, it follows that $(L^\infty_x)_{x\in(0,M]}$ has the same 
distribution as $(\Lambda^*(\tau^*_x))_{x\in(0,M]}$, giving the first assertion of the proposition. The second
assertion follows by differentiating the mapping $x\mapsto \Lambda^*(\tau^*_x)$.  \endproof

\section{The stochastic differential equation for areas of spheres}
\label{sec:EDS}

With a slight abuse of notation, for every
 $w\in \R$ and $\lambda>0$, we write $\P_{w,\lambda}$ for a probability measure under which 
the pair $(Z_s,\Lambda_s)_{s\geq 0}$ is such that
$(Z_s)_{s\geq 0}$ is distributed as the solution of \eqref{EDS} with initial value $w$, and
$$\Lambda_s=\lambda\, \exp\Big(-\int_0^s Z_r\,\dd r\Big).$$

\begin{proposition}
\label{transker}
The process $(W^*(\tau^*_x),\Lambda^*(\tau^*_x))_{x>0}$ is a time-homogeneous Markov process, whose
transition kernels $(\Pi_s)_{s>0}$ are given by
$$\int \Pi_s((w,\lambda), \dd(w',\lambda'))\,\varphi(w',\lambda')= \E_{w,\lambda}[\varphi(Z_{\eta_s},\Lambda_{\eta_s})],$$
where
$$\eta_s:=\inf\Big\{r\geq 0: \int_0^r (\Lambda_u)^{1/3}\,\dd u \geq s\Big\}.$$
\end{proposition}

\proof We fix $a>0$ and use the notation
$$T^{(a)}:=\sup\Big\{t\geq 0: \int_0^t W^*(u)\,\dd u\geq -a\Big\}.$$
Note in particular that 
$\int_0^{T^{(a)}} W^*(u)\,\dd u=-a$. By Proposition \ref{tecpro}, we have
$$\Big(W^*(T^{(a)}-s)\Big)_{0\leq s\leq T^{(a)}}\build{=}_{}^{\rm (d)} \Big(W^*(s)\Big)_{0\leq s\leq T^{(a)}}.$$
We can in fact extend this identity as 
\begin{equation}
\label{reverstar}
\Big(W^*(T^{(a)}-s)\Big)_{s\geq 0}\build{=}_{}^{\rm (d)} \Big(W^*(s)\Big)_{s\geq 0}.
\end{equation}
To justify this, let $b>0$, and recall the definition \eqref{defSb} of $S_b$. Using the
fact that $(W^*(S_b+t))_{t\geq0}$ has the same distribution as $(W^*(s))_{s\geq 0}$, 
and Proposition \ref{tecpro} again, we get
$$\Big(W^*(T^{(a)}-s)\Big)_{0\leq s\leq T^{(a)}-S_b} \build{=}_{}^{\rm (d)} \Big(W^*(T^{(a+b)}-s)\Big)_{0\leq s\leq T^{(a+b)}}\
\build{=}_{}^{\rm (d)} \Big(W^*(s)\Big)_{0\leq s\leq T^{(a+b)}}.$$
(For the first identity, note that, if $\wt W(u):=W^*(S_b+u)$, we have $\wt T^{(a+b)}=T^{(a)}-S_b$, with an obvious
notation, and $\wt W(\wt T^{(a+b)}-s)=W^*(T^{(a)}-s)$.) Then we just have to let $b\uparrow \infty$
to get \eqref{reverstar}. 

We then set, for every $s\geq 0$,
\begin{align}
\label{tec66}
\ov \Lambda^{(a)}(s):= \Lambda^*(T^{(a)}-s)&=\exp\Big(\int_0^{T^{(a)}-s} W^*(r)\,\dd r\Big)\nonumber\\
&=e^{-a}\exp\Big(\int_{T^{(a)}}^{T^{(a)}-s} W^*(r)\,\dd r\Big)\nonumber\\
&=e^{-a}\exp\Big(-\int_{0}^{s} W^*(T^{(a)}-r)\,\dd r\Big)\nonumber\\
&=e^{-a}\exp\Big(-\int_{0}^{s} \ov Z^{(a)}(r)\,\dd r\Big)
\end{align}
where we have set $\ov Z^{(a)}(r):= W^*(T^{(a)}-r)$ for every $r\geq 0$. 

From \eqref{reverstar} and \eqref{tec66}, we
see that $(\ov Z^{(a)}(r),\ov \Lambda^{(a)}(r))_{r\geq 0}$
has the same distribution as $(Z_r,\Lambda_r)_{r\geq 0}$ under
$\int\mu(\dd w)\,\P_{w,e^{-a}}$. We then set
$$I^{(a)}:=\int_{T^{(a)}}^\infty  \Lambda^*(s)^{1/3}\,\dd s.$$
Let $x>0$. On the event where $I^{(a)}\leq x$, we have
$$\tau^*_x=\sup\Big\{z\leq T^{(a)}:\int_z^\infty \Lambda^*(s)^{1/3}\,\dd s\geq x\Big\}$$
and thus
\begin{align*}
T^{(a)}-\tau^*_x&=\inf\Big\{s\geq 0: \int_{T^{(a)}-s}^{T^{(a)}} \Lambda^*(u)^{1/3}\,\dd u\geq x - \int_{T^{(a)}}^{\infty} \Lambda^*(u)^{1/3}\,\dd u\Big\}\\
&=\inf\Big\{s\geq 0: \int_0^s \ov\Lambda^{(a)}(u)^{1/3}\,\dd u\geq x-I^{(a)}\Big\}.
\end{align*}
Still on the same event, we have
$$\Lambda^*(\tau^*_x)= \ov\Lambda^{(a)}(T^{(a)}-\tau^*_x)
=\ov\Lambda^{(a)}\Big(\inf\Big\{s\geq 0: \int_0^s \ov\Lambda^{(a)}(u)^{1/3}\,\dd u\geq x-I^{(a)}\Big\}\Big)$$
and, if we write $\eta^{(a)}_x=T^{(a)}-\tau^*_x$, we also get, for every $y>x$,
\begin{align}
\label{tec67}
\Lambda^*(\tau^*_y)=\ov\Lambda^{(a)}(\eta^{(a)}_y)&=\ov\Lambda^{(a)}\Big(\inf\Big\{s\geq 0: \int_0^s \ov\Lambda^{(a)}(u)^{1/3}\,\dd u\geq y-I^{(a)}\Big\}\Big)
\nonumber\\
&=\ov\Lambda^{(a)}\Big(\eta^{(a)}_x+\inf\Big\{s\geq 0: \int_{\eta^{(a)}_x}^{\eta^{(a)}_x+s} \ov\Lambda^{(a)}(u)^{1/3}\,\dd u\geq y-x\Big\}\Big)\nonumber\\
&=\ov\Lambda^{(a)}\Big(\eta^{(a)}_x+\inf\Big\{s\geq 0: \int_0^{s} \ov\Lambda^{(a)}(\eta^{(a)}_x+u)^{1/3}\,\dd u\geq y-x\Big\}\Big),
\end{align}
and similarly
\begin{equation}
\label{tec68}
W^*(\tau^*_y)=\ov Z^{(a)}(\eta^{(a)}_y)=\ov Z^{(a)}\Big(\eta^{(a)}_x+\inf\Big\{s\geq 0: \int_0^{s} \ov\Lambda^{(a)}(\eta^{(a)}_x+u)^{1/3}\,\dd u\geq y-x\Big\}\Big).
\end{equation}
Let us then consider the filtration $(\g_t)_{t\geq 0}$ defined by
$$\g_t:=\sigma\Big(\ov Z^{(a)}(r),\ov\Lambda^{(a)}(r):0\leq r\leq t\Big)\vee \sigma\Big(\Lambda^*(T^{(a)}+s):s\geq 0\Big).$$
By convention, we set $\eta^{(a)}_x=\infty$ if $I^{(a)}>x$. Note that the event $\{I^{(a)}\leq x\}$ is $\g_0$-measurable.

By \eqref{reverstar}, $(\ov Z^{(a)}(t))_{t\geq 0}$ has the same distribution as $(W^*(t))_{t\geq 0}$ and 
is thus distributed as the solution of \eqref{EDS} with initial distribution $\mu$. It follows that $(\ov Z^{(a)}(t),\ov\Lambda^{(a)}(t))_{t\geq 0}$ is 
(time-homogeneous) Markov with respect to the filtration $(\g_t)_{t\geq 0}$, and its transition kernels $\Theta_s$, $s\geq 0$,
are specified by saying that $\Theta_s((w,\lambda),\cdot)$ is the distribution of $(Z_s,\Lambda_s)$ under $\P_{w,\lambda}$.
Moreover $(\ov Z^{(a)}(0),\ov\Lambda^{(a)}(0))$ is distributed according to $\mu(\dd w)\delta_{e^{-a}}(\dd\lambda)$, by \eqref{tec66}. 

Observe that $\eta^{(a)}_x$ is a stopping time of the filtration $(\g_t)_{t\geq 0}$, and that, on the
event where $I^{(a)}\leq x$, \eqref{tec67} and \eqref{tec68} express $(\ov Z^{(a)}(\eta^{(a)}_y),\ov\Lambda^{(a)}(\eta^{(a)}_y))$
in terms of $\eta^{(a)}_x$ and the shifted process 
$$(\ov Z^{(a)}(\eta^{(a)}_x+u),\ov\Lambda^{(a)}(\eta^{(a)}_x+u))_{u\geq 0}.$$
 An application 
of the strong Markov property shows that, on the event $\{I^{(a)}\leq x\}$, conditionally on the 
$\sigma$-field $\g_{\eta^{(a)}_x}$, the distribution of the pair $(\ov Z^{(a)}(\eta^{(a)}_y),\ov\Lambda^{(a)}(\eta^{(a)}_y))$
is the law of $(Z_{\eta_{y-x}},\Lambda_{\eta_{y-x}})$ under $\P_{\ov Z^{(a)}(\eta^{(a)}_x),\ov\Lambda^{(a)}(\eta^{(a)}_x)}$,
which by definition is
$$\Pi_{y-x}\Big((\ov Z^{(a)}(\eta^{(a)}_x),\ov\Lambda^{(a)}(\eta^{(a)}_x)),\,\cdot\,\Big).$$

Consider now $0<x_1<x_2<\cdots<x_k=x$ and $y>x$. By the preceding observations, we have
\begin{align*}
&\E\Big[ \varphi\Big(W^*(\tau^*_y), \Lambda^*(\tau^*_y)\Big)\,F\Big((W^*(\tau^*_{x_j}), \Lambda^*(\tau^*_{x_j}))_{1\leq j\leq k}\Big)\,\mathbf{1}_{\{I^{(a)}\leq x_1\}}\Big]\\
&\qquad=\E\Big[ \varphi\Big(\ov Z^{(a)}(\eta^{(a)}_y),\ov \Lambda^{(a)}(\eta^{(a)}_y)\Big)\,F\Big((\ov Z^{(a)}(\eta^{(a)}_{x_j}), \ov \Lambda^{(a)}(\eta^{(a)}_{x_j}))_{1\leq j\leq k}\Big)\,\mathbf{1}_{\{I^{(a)}\leq x_1\}}\Big]\\
&\qquad=\E\Big[ \Pi_{y-x}\varphi\Big(\ov Z^{(a)}(\eta^{(a)}_x),\ov \Lambda^{(a)}(\eta^{(a)}_x)\Big)\,F\Big((\ov Z^{(a)}(\eta^{(a)}_{x_j}), \ov \Lambda^{(a)}(\eta^{(a)}_{x_j}))_{1\leq j\leq k}\Big)\,\mathbf{1}_{\{I^{(a)}\leq x_1\}}\Big]\\
&\qquad=\E\Big[ \Pi_{y-x}\varphi\Big(W^*(\tau^*_x), \Lambda^*(\tau^*_x)\Big)\,F\Big((W^*(\tau^*_{x_j}), \Lambda^*(\tau^*_{x_j}))_{1\leq j\leq k}\Big)\,\mathbf{1}_{\{I^{(a)}\leq x_1\}}\Big]
\end{align*}
and we just have to let $a\uparrow \infty$ to get the Markov property as stated in the proposition. \endproof

We can now easily derive Theorem \ref{Markov-LT}.

\proof[Proof of Theorem \ref{Markov-LT}] 
Recall that $L^\infty_x>0$, for every $x>0$, a.s. From Propositions \ref{repre-loc} and \ref{transker}, we readily obtain
that $(L^\infty_x,\dot L^\infty_x)_{x> 0}$ is Markov with transition kernels given by
$$\mathbf{\Pi}_s((\ell,z),A)= \int \Pi_s((-\ell^{-2/3}z,\ell),\dd(w',\lambda'))\,\mathbf{1}_A(\lambda',-\lambda'^{2/3}w'),$$
for every $\ell >0$, $z\in\R$ and every Borel subset $A$ of $(0,\infty)\times \R$. Since we have also 
$L^\infty_0=\dot L^\infty_0=0$, this suffices to get the desired Markov property. \endproof

Recall the notation $(\f_t)_{t\geq 0}$ for the (completion of the) canonical filtration of $(L^\infty_t)_{t\geq 0}$. For $\ve>0$, we 
also set $\f^\ve_t:=\f_{\ve+t}$ for every $t\geq 0$. 

\begin{proposition}
\label{main-theo}
Let $\ve>0$. Then the process $(L^\infty_{\ve+t},\dot L^\infty_{\ve+t})_{t\geq 0}$ 
satisfies the stochastic differential equation
$$\dot L^\infty_{\ve+t}=\dot L^\infty_\ve + 4\int_0^t \sqrt{L^\infty_{\ve+s}}\,\dd B^\ve_s + \int_0^t h(L^\infty_{\ve+s},\dot L^\infty_{\ve+s})\,\dd s,$$
where $B^\ve$ is an $(\f^\ve_t)$-Brownian motion started at $0$, and
$$h(t,y)= -8t\,\frac{p'_t(-y/2)}{p_t(-y/2)} + \frac{4}{3}\,\frac{y^2}{t}.$$
\end{proposition}

\proof If $w\in\R$ and $\lambda>0$, the process $(Z_{\eta_s},\Lambda_{\eta_s})_{s\geq 0}$ is under $\P_{w,\lambda}$ time-homogeneous
Markov with the transition kernels $\Pi_s$ of Proposition \ref{transker}. We can thus combine 
Propositions \ref{repre-loc} and \ref{transker} to get that $(L^\infty_{\ve+t},\dot L^\infty_{\ve+t})_{t\geq 0}$  has the same 
distribution as 
$$\Big(\Lambda_{\eta_s},-(\Lambda_{\eta_s})^{2/3}Z_{\eta_s}\Big)_{s\geq 0}$$
under $\P_{W^*(\tau^*_\ve),\Lambda^*(\tau^*_\ve)}$. By construction, under the probability measure
$\P_{W^*(\tau^*_\ve),\Lambda^*(\tau^*_\ve)}$, the pair $(Z_s,\Lambda_s)$ satisfies the equation
\begin{align*}
\dd Z_s&=4\,\dd B_s + b(Z_s)\,\dd s\\
\dd \Lambda_s&=-Z_s\Lambda_s\,\dd s.
\end{align*}
If $X_s:=-(\Lambda_s)^{2/3}Z_s$, an application of It\^o's formula gives
\begin{align*}
\dd X_s&=-4(\Lambda_s)^{2/3}\,\dd B_s - (\Lambda_s)^{2/3}\,b(Z_s)\,\dd s + \frac{2}{3}\, (\Lambda_s)^{-1/3}\,(Z_s)^2\,\Lambda_s\,\dd s\\
&=-4\,(\Lambda_s)^{2/3}\,\dd B_s + (\Lambda_s)^{1/3}\,h(\Lambda_s,X_s)\,\dd s,
\end{align*}
where we used the scaling relations
$$p_t(z)=t^{-2/3}p_1(t^{-2/3}z)\;,\quad p'_t(z)=t^{-4/3}p'_1(t^{-2/3}z)$$
to verify that, for every $z\in\R$ and $t>0$,
$$-t^{1/3}b(z)+\frac{2}{3}\, t^{1/3} z^2=-8t^{1/3}\,\frac{p'_1(z/2)}{p_1(z/2)}\,+\frac{4}{3}\,t^{1/3}z^2=h(t,-t^{2/3}z).$$

We then use a standard time change argument. If 
$$M_t:= -\int_0^t (\Lambda_s)^{1/6}\,\dd B_s,$$
the Dubins-Schwarz theorem \cite[Theorem V.1.6]{RY} shows that $\beta_t:=M_{\eta_t}$ is a linear Brownian motion (in the appropriate time-changed filtration). Moreover,
$$-\int_0^{\eta_t} (\Lambda_s)^{2/3}\,\dd B_s =\int_0^{\eta_t}\sqrt{\Lambda_s}\,\dd M_s= \int_0^t\sqrt{\Lambda_{\eta_u}}\,\dd \beta_u,$$
where the last equality follows from the formal change of variables $s=\eta_u$, which is easily justified (we omit the details). Since we have also
$$\int_0^{\eta_t} (\Lambda_s)^{1/3} \,h(\Lambda_s,X_s)\,\dd s= \int_0^t h(\Lambda_{\eta_u},X_{\eta_u})\,\dd u,$$
we conclude that
\begin{equation}
\label{EDS29}
X_{\eta_t}= X_0+ 4 \int_0^t\sqrt{\Lambda_{\eta_u}}\,\dd \beta_u + \int_0^t h(\Lambda_{\eta_u},X_{\eta_u})\,\dd u.
\end{equation}

It is now easy to see that the pair $(L^\infty_{\ve+t},\dot L^\infty_{\ve+t})_{t\geq 0}$
satisfies the equation stated in the theorem. To justify this, recall that $\E[L^\infty_x]=c\,x^3$ with a constant $c<\infty$,
Since $(L^\infty_{\ve+t},\dot L^\infty_{\ve+t})_{t\geq 0}$  has the same distribution as $(\Lambda_{\eta_t},X_{\eta_t})_{t\geq 0}$, 
it follows from \eqref{EDS29} that
$$\dot L^\infty_{\ve+t}-\dot L^\infty_{\ve}-\int_0^t h(L^\infty_{\ve+s},\dot L^\infty_{\ve+s})\,\dd s$$
is a continuous martingale in the filtration $(\f^\ve_t)_{t\geq 0}$, with quadratic variation 
$$16\,\int_0^t L^\infty_{\ve+s}\,\dd s.$$
This martingale can be represented as 
$$4\,\int_0^t \sqrt{L^\infty_{\ve+s}}\,\dd B^\ve_s,$$
where $B^\ve$ is an $(\f^\ve_t)$-Brownian motion started from $0$. 
 This completes the proof. \endproof
 
 We finally prove Theorem \ref{mainT}, which was stated in the introduction. 

\proof[Proof of Theorem \ref{mainT}]
We first observe that $h(t,x)>0$ for every $t>0$ and $x\in\R$. As it is clear from \eqref{formupt} and \eqref{map-Airy},
the function $x\mapsto  x\mathrm{Ai}(x^2) + \mathrm{Ai}'(x^2)$ is everywhere negative, and since 
the function $x\mapsto \mathrm{Ai}(x^2)$ is positive on $\R$, it follows from formula \eqref{formu-drift} that $h(t,x)>0$.

Recall from Proposition \ref{localtimesBP} that $L^\infty_t>0$ for every $t>0$, a.s. 
Then, for $\ve>0$, the Brownian motion $B^\ve$ in Proposition \ref{main-theo}
can be written as
$$B^\ve_t=\int_0^t \frac{1}{4\sqrt{L^\infty_{\ve+s}}}\,\dd H^\ve_s,$$
where
$$H^\ve_t=\dot L^\infty_{\ve+t} - \int_0^t h(L^\infty_{\ve+s},\dot L^\infty_{\ve+s})\,\dd s.$$
If $0<\ve'<\ve$, it is a simple matter to verify that we must have $B^\ve_t=B^{\ve'}_{(\ve-\ve')+t}-B^{\ve'}_{\ve-\ve'}$,
and then that there exists an $(\f_t)$-Brownian motion $B$ started at $0$
such that, for every $\ve>0$, $B^\ve_t=B_{\ve+t}-B_\ve$. 
From Proposition \ref{main-theo}, we then get, for every fixed 
$u>0$ and $0<\ve<u$,
\begin{equation}
\label{EDSeps}
\dot L^\infty_u =\dot L^\infty_\ve + 4\int_\ve^u \sqrt{L^\infty_s}\,\dd B_s +\int_\ve^u h(L^\infty_s,\dot L^\infty_s)\,\dd s.
\end{equation}
By Proposition \ref{localtimesBP}, $\dot L^{\infty}_\ve$ tends to $0$ as $\ve\to 0$, and thus we have
$$\dot L^\infty_u -\dot L^\infty_\ve - 4\int_\ve^u \sqrt{L^\infty_s}\,\dd B_s \la \dot L^\infty_u - 4\int_0^u \sqrt{L^\infty_s}\,\dd B_s$$
in probability, as $\ve \to 0$.  Hence, 
$$\int_\ve^u h(L^\infty_s,\dot L^\infty_s)\,\dd s$$
also converges in probability to a finite limit when $\ve\to 0$. This implies that
$$\int_0^u h(L^\infty_s,\dot L^\infty_s)\,\dd s <\infty,\;\hbox{a.s.}$$
and, to complete the proof of \eqref{EQ1}, we just have to let $\ve\to 0$ in \eqref{EDSeps}.
\endproof

\noindent{\bf Acknowledgments.} I am indebted to an anonymous referee whose numerous remarks allowed me
to improve the original manuscript.
 I also thank Patrick Cattiaux for useful comments about the time reversal of diffusions.

\end{document}